\documentclass{amsart}
\usepackage{amsfonts, amssymb, latexsym, amscd}
\usepackage{amsmath}
\usepackage[mathscr]{eucal}
\usepackage{enumerate}
\usepackage[headings]{fullpage}
\usepackage[pdftex]{graphicx, color}
\usepackage[usenames,dvipsnames]{xcolor}
\theoremstyle{plain}
\newtheorem{theorem}{Theorem}[section]
\newtheorem{corollary}[theorem]{Corollary}
\newtheorem{lemma}[theorem]{Lemma}
\newtheorem{proposition}[theorem]{Proposition}

\theoremstyle{definition}

\newtheorem{remark}[theorem]{Remark}

\newcommand{\pref}[1]{(\ref{#1})}

\newcommand{\RR}{\mathbb{R}}

\newcommand{\MM}{\mathcal{M}}

\newcommand{\vp}{\varphi}

\newcommand{\al}{\alpha}

\newcommand{\ga}{\gamma}

\newcommand{\ka}{\kappa}
\newcommand{\la}{\lambda}

\newcommand{\thet}{\theta}

\newcommand{\Ga}{\Gamma}

\newcommand{\La}{\Lambda}

\newcommand{\norm}[1]{\left|\left|#1\right|\right|}
 
\newcommand{\rst}[1]{\ensuremath{{\mathbin\upharpoonright}%
\raise-.5ex\hbox{$#1$}}}

\newcounter{gscan}
\newcounter{btscan}
\newcounter{cscan}
\newcounter{hscan}
\newcounter{bhscan}

\numberwithin{equation}{section}

\begin{document}

\title[A Quantitative Regularity Estimate for Parabolic Equations]{A Quantitative Regularity Estimate for Nonnegative Supersolutions of Fully Nonlinear Uniformly Parabolic Equations}

\author{Jessica Lin}
\address{University of Chicago\\
Department of Mathematics\\
Chicago, IL 60637}
\email[Jessica Lin]{jessica@math.uchicago.edu}
\urladdr{http://math.uchicago.edu/$\sim$jessica}
\date{\today}

\begin{abstract}
This note establishes an interior quantitative lower bound for nonnegative supersolutions of fully nonlinear uniformly parabolic equations. The result may be interpreted as a nonlinear, quantitative version of a growth lemma established by Krylov and Safonov for nonnegative supersolutions of linear uniformly parabolic equations in nondivergence form. Our approach is different, and follows from an application of the Fabes-Stroock estimate. The result is the parabolic analogue of an elliptic regularity estimate established by Caffarelli, Souganidis, and Wang in the stochastic homogenization of fully nonlinear uniformly elliptic equations. 
\end{abstract}

\maketitle

\section{Introduction}
An interesting question in the theory of elliptic and parabolic partial differential equations is whether it is possible to have a comparison between a function $u$ and $\mathcal{L}u$, where $\mathcal{L}$ is a fully nonlinear uniformly elliptic differential operator. In the time-independent setting of $\RR^{N}$, this question has been completely addressed. The Alexandroff-Backelman-Pucci estimate \cite{cc} yields that if $\mathcal{L}u\geq f$ in a cylinder and $u=0$ on the parabolic boundary, then $u$ is controlled from above by the $L^{N}$-norm of $f$. Moreover, it was shown in \cite{csw} that under some additional assumptions on $u$ and $f$, $u$ in the interior is controlled from below by the $L^{N}$-norm of $f$ raised to a power depending only on the ellipticity constants and dimension of the space. In the context of uniformly parabolic equations in $\RR^{N+1}$, the parabolic analogue of the Alexandroff-Backelman-Pucci estimate was first established by Krylov \cite{krylovabp} for $u\in W^{2,1,N+1}$ and revisited by Tso \cite{tso} with a simplified approach of proof. The parabolic estimate is now referred to as the Alexandroff-Backelman-Pucci-Krylov-Tso estimate, and it was extended to viscosity solutions by Wang in \cite{wangreg1}. The purpose of this paper is to address how $u$ is controlled from below by the $L^{N+1}$-norm of $f$. 

We prove here that in a fraction of the original domain, nonnegative supersolutions of uniformly parabolic equations are bounded below by $\norm{f}_{L^{N+1}}^{\al}$, where $\al\sim\frac{1}{\norm{f}^{2(N+1)}_{L^{N+1}}\left|\log (\norm{f}_{L^{N+1}})\right|}$. Although the result is presented for fully nonlinear equations, it is to our knowledge new for linear nondivergence form equations with bounded measurable coefficients. Our initial motivation was to develop a parabolic version of the lower bound established by Caffarelli, Souganidis, and Wang in \cite{csw}, which was used in the error estimates for stochastic homogenization of uniformly elliptic equations in random media \cite{cs}. Although our general approach follows theirs, it is necessary to develop a number of new arguments to handle the parabolic structure of the problem. We also show that we recover the elliptic result of \cite{csw} from our estimates in the limit as $t\rightarrow\infty$. 

Before stating the result, we briefly explain the notation and setting. We consider $u\in C(\RR^{N+1})$ satisfying in the viscosity sense
\begin{equation}\label{nonlinear}
\begin{cases}
u_{t}-F(D^{2}u, x, t)\geq f\quad\text{in}\quad Q_{1},\\
u=0\quad\text{on}\quad \partial_{p}Q_{1}\\
u\geq 0
\end{cases}
\end{equation}
where $Q_{1}$ and $\partial_{p}Q_{1}$  represent the unit cylinder and its parabolic boundary. That is, 
\begin{align*}
Q_{1}&=B_{1}(0)\times (-1, 0]\subset \RR^{N+1}, \quad \text{and} \\
\partial_{p}Q_{1}&=\left(B_{1}(0)\times \left\{t=-1\right\}\right)\cup \left(\partial B_{1}(0)\times [-1, 0)\right).
\end{align*}
In general, we use the notation
\begin{align*}
Q_{r}(x_{0}, t_{0})&=B_{1}(x_{0})\times (t_{0}-r^{2}, t_{0}]\subset \RR^{N+1},\quad \text{and}\\
\partial_{p}Q_{r}(x_{0}, t_{0})&=\left(B_{r}(x_{0})\times \left\{t=t_{0}-r^{2}\right\}\right)\cup \left(\partial B_{r}(x_{0})\times [t_{0}-r^{2}, t_{0})\right).
\end{align*}
We will frequently refer to $u$ satisfying \pref{nonlinear} as a supersolution to $u_{t}-F(D^{2}u, x, t)=f$. 

We assume that $F$ is uniformly elliptic, with ellipticity constants $\la, \La$, that is, for every $M,K\in \mathbb{S}^{N}$(the space of $N\times N$ symmetric matrices), $K\geq 0$, we have 
\begin{equation}\label{ellipcond}
\la \norm{K}\leq F(M+K, x, t)-F(M, x,t)\leq \La\norm{K}\quad\text{for all}\quad (x,t)\in Q_{1},
\end{equation}
where $\norm{K}$ denotes the maximum eigenvalue of $K$. This is equivalent to saying that \pref{nonlinear} is uniformly parabolic. In addition, we assume that $F$ satisfies the standard regularity assumptions to admit a comparison principle (see \cite{users, lpvisc}).

The main result is:
\begin{theorem}\label{thmnonlin}
Fix $f$ so that $0\leq f+F(0,\cdot, \cdot)\leq 1$ and assume \pref{ellipcond}. Let $u$ satisfy \pref{nonlinear}. For every $\ka\in(0,1)$, there exists $c, C, \rho, \beta> 0$, depending only on $\la, \La, N, \ka$, such that for all $|x|\leq \ka$ and $0\geq t\geq -\frac{\ka}{2|Q_{1}|}\norm{f+F(0,x,t)}^{N+1}_{L^{N+1}(Q_{1})}$,
\begin{align*} 
&c\norm{f+F(0, x,t)}^{\rho}_{L^{N+1}(Q_{1})}\exp (-\beta\norm{f+F(0,x,t)}_{L^{N+1}(Q_{1})}^{-2(N+1)}) \leq u(x,t)\\
&\leq C \norm{f+F(0,x,t)}_{L^{N+1}(Q_{1})}.
\end{align*}
\end{theorem}


As previously mentioned, the upper bound is the Alexandroff-Backelman-Pucci-Krylov-Tso estimate \cite{wangreg1}. The focus of this note is to obtain the lower bound. We appeal to the standard ``linearized" interpretation of \pref{nonlinear} using Pucci's extremal operators (see \cite{cc, wangreg1}). For every $M\in \mathbb{S}^{N}$, the lower and upper Pucci's extremal operators are
\begin{equation*}
\MM^{-}(M, \la, \La):=\MM^{-}(M)=\la\left(\sum_{e_{i}>0} e_{i}\right)+\La\left(\sum_{e_{i}<0}e_{i}\right),
\end{equation*}
and
\begin{equation*}
\MM^{+}(M, \la, \La):=\MM^{+}(M)=\La\left(\sum_{e_{i}>0} e_{i}\right)+\la\left(\sum_{e_{i}<0}e_{i}\right),
\end{equation*}
where $e_{i}$ are the eigenvalues of $M$. It follows that $\MM^{-}(\cdot)$ and $\MM^{+}(\cdot)$ are both uniformly elliptic. 
For more properties of Pucci's extremal operators, see \cite{cc, wangreg1}. 

 
For any domain $D\subset \RR^{N+1}$, we define $\overline{S}(f, D)$ to be the collection of $u\in C(\RR^{N+1})$ satisfying in the viscosity sense,
\begin{equation*}
u_{t}-\MM^{-}(D^{2}u)\geq f\quad\text{in}\quad D,
\end{equation*}
Respectively, we define $\underline{S}(f, D)$ to be the collection of $u\in C(\RR^{N+1})$ satisfying in the viscosity sense,
\begin{equation*}
u_{t}-\MM^{+}(D^{2}u)\leq f\quad\text{in}\quad D.
\end{equation*}
Although the original formulation introduced in \cite{cc} is for $f\in C(\RR^{N+1})$, the theory of viscosity solutions for $f\in L^{\infty}(\RR^{N+1})$ and compactly supported is established in the work of \cite{lpvisc}, and the properties in \cite{cc, wangreg1} easily generalize to this setting. 

It is shown in \cite[Lemma 3.12]{wangreg1}  that $u$ satisfying \pref{nonlinear} also satisfies $u\in \overline{S}(f+F(0,\cdot,\cdot), Q_{1})$. Therefore, in order to prove Theorem \ref{thmnonlin}, it is enough to prove
\begin{theorem} \label{TSFS}
Fix $0\leq f\leq 1$, and let $u$ be nonnegative in $\overline{S}(f, Q_{1})$, with $u=0$ on $\partial_{p}Q_{1}$. For every $\ka\in(0,1)$, there exists $c, C, \rho, \beta>0$ which depend only on $\la, \La, N, \ka$, so that for all $|x|\leq \ka$ and $0\geq t \geq -\frac{\ka}{2|Q_{1}|}\norm{f}^{N+1}_{L^{N+1}(Q_{1})}$,
\begin{equation}\label{bothsides}
c\norm{f}_{L^{N+1}(Q_{1})}^{\rho}\exp(-\beta \norm{f}^{-2(N+1)}_{L^{N+1}(Q_{1})})\leq u(x,t)\leq C \norm{f}_{L^{N+1}(Q_{1})}.
\end{equation}
\end{theorem}


 


We note that in the parabolic setting, the domain where \pref{bothsides} holds depends on $\norm{f}_{L^{N+1}(Q_{1})}$. This is a consequence of the causality property inherent to solutions of parabolic equations. If $u$ solves a parabolic equation, $u(\cdot,t)$ is only affected by $f(\cdot, s)$ for $s\leq t$. Thus, any estimate for $u$ in terms of $f$ will hold for times after we ``see" $f$. In general, the domain where $f$ is large is comparable to $\frac{1}{|Q_{1}|}\norm{f}^{N+1}_{L^{N+1}(Q_{1})}$, and this explains the dependence in \pref{bothsides}. In Section \ref{secext}, we present several special cases of Theorem  \ref{TSFS} which yield sharper estimates given more information about the distribution of $f$. 


Theorem \ref{TSFS} follows relatively easily from:
\begin{theorem}\label{LB}
Let $u$ be nonnegative in $\overline{S}(f, Q_{1})$. Set $\Ga=\left\{f>\al\right\}\subset Q_{1}$, and $m=\frac{|\Ga|}{|Q_{1}|}$, where $|\cdot|$ denotes the Lebesgue measure on $\RR^{N+1}$. For every $\ka\in(0,1)$, there exists $c, \rho, \beta>0$ depending only $\la, \La, N, \ka$, so that for all $|x|\leq \ka$ and $0\geq t\geq -\ka m$,
\begin{equation}\label{measlb}
u(x,t)\geq cm^{\rho}\exp(-\beta/m^{2})\al.
\end{equation}
\end{theorem}

In order to prove Theorem \ref{LB}, we compare $u\in \overline{S}(f, Q_{1})$ to $w(x,t)=w(x,t; \Ga)$, the ``fundamental solution corresponding to the domain $\Ga$," which solves
\begin{equation}\label{nonlinearfund}
\begin{cases}
w_{t}-\MM^{-}(D^{2}w)=\chi_{\Ga}\quad\text{in}\quad Q_{1},\\
w=0\quad\text{on}\quad\partial_{p}Q_{1},
\end{cases}
\end{equation}
where $\chi_{\Ga}$ denotes the characteristic function of a measurable set $\Ga\subset Q_{1}$. We first prove \pref{measlb} for the solution $w$. Theorem \ref{LB} may be interpreted as a quantitative version of a growth lemma established by Krylov and Safonov (see \cite{krylov}, Theorem 4.2.1 and \cite{krylovsaf2},\S 8), where the lower bound is given in terms of an unknown function $\vp(|\Ga|)$ satisfying $\vp(m)>0$ for $m>0$. Our approach differs from the proof presented in \cite{krylov}, where the author uses classical covering arguments to cover $\Ga$ with cylinders which contain a significant proportion of $\Ga$ in measure. A difficulty in quantifying this argument comes from the fact that these cylinders often spill outside of the original domain $Q_{1}$. Under an additional assumption that $\Ga\subset Q_{r}(0,-1/2)$ for some $r<1$, this difficulty is avoided and quantitative estimates have been studied in this setting.  In \cite{fabgarsalsa}, the authors obtain a quantitative lower bound for the Green's function, which yields an estimate for solutions of \pref{nonlinearfund}. In \cite{krylovL}, the author presents a quantitative lower bound for solutions of \pref{nonlinear} directly, without appealing the Green's function representation. Under this additional assumption on $\Ga$, the estimates in \cite{fabgarsalsa} and \cite{krylovL} are stronger than those in Proposition \ref{qlbnd} and Theorem \ref{LB}, respectively. In Section \ref{secext}, we recover these stronger estimates as a corollary to our results.

Our approach is independent of the works previously mentioned. We construct a covering  of $\Ga$ which is completely contained inside of $Q_{1}$, and this allows us to obtain a quantitative lower bound under the more general hypotheses first used in \cite{krylov} and \cite{krylovsaf2}. The Fabes-Stroock estimate \cite{fs, amarnorando} is applied to compare fundamental solutions corresponding to different domains. In particular, we control $w(x, t; \Ga\cap Q_{r})$ from below by $w(x,t; Q_{r})$ for some $Q_{r}$ which contains a significant proportion of $\Ga$ in measure. In order to control $w(x,t; Q_{r})$ from below, we introduce an iterative method to prove regularity estimates in larger space domains for later times. Moreover, our approach can be adapted to special cases which have various applications. Corollary \ref{homogest} in particular is the key regularity estimate utilized in the study of error estimates for stochastic homogenization of uniformly parabolic equations  \cite{linhomog1}. 


This paper is organized as follows. Section \ref{seccyl} is devoted to establishing a quantitative lower bound for $w(x,t; Q_{r})$. We revisit some of the results in various works by Krylov \cite{krylov} and Krylov and Safonov \cite{krylovsaf, krylovsaf2}, relaxing some of the hypotheses and presenting the proofs for Pucci's extremal operators. We also describe the iterative construction, which is completely contained inside of $Q_{1}$ and allows us to obtain regularity estimates in larger space domains at later times. In Section \ref{secgen}, we use a covering argument to complete the proof of Theorem \ref{LB}. We consequently obtain Theorem \ref{TSFS}. In Section \ref{secext}, we present several special cases of Theorem \ref{LB}, which follow easily by modifying our approach with the new information provided. 





\section{Quantitative Lowers bounds on Fundamental Solutions for Subcylinders}\label{seccyl}

We prove an interior quantitative lower bound for fundamental solutions of subcylinders, $w(x,t; Q_{r}(x_{0}, t_{0}))$, with $Q_{r}(x_{0}, t_{0})\subset Q_{1}$. We note that in the elliptic setting, one can generally obtain an interior lower bound for $w(x; B_{r}(x_{0}))$  by iteratively applying the weak Harnack inequality. However, in the parabolic setting, the argument becomes more delicate in order to account for the time shifts of the parabolic Harnack inequality \cite{wangreg1}. Our proof is inspired by the iterative approach in the elliptic setting, however we do not employ the parabolic Harnack inequality. Instead, we use some of the ideas and constructions originally found in \cite{krylov} to compare supersolutions in towers of oblique cylinders. Our construction utilizes lower bounds from previous times to obtain lower bounds at later times, and manages to stay inside of $Q_{1}$ the duration of the process. We point out that this section is completely self-contained, and does not require any additional knowledge of parabolic regularity theory.

We first present an important comparison lemma found in \cite{krylovsaf, krylov}, which allows us to obtain a ``weak" comparison principle for oblique cylinders. Although the original estimate in \cite{krylovsaf, krylov} holds for more general linear operators of the form $D_{t}-a_{ij}(x,t)\frac{\partial^{2}}{\partial x_{i} \partial x_{j}}-b_{i}(x,t)\frac{\partial}{\partial x_{i}}+c(x,t)$, with $\left\{a_{ij}\right\}$ uniformly elliptic, $\left\{b_{i}\right\}$ uniformly bounded, $c$ bounded and nonnegative, the original estimate also requires stronger hypotheses on the dimensions of the oblique cylinder. Our presentation here relaxes these hypotheses and gives more detailed estimates on how the constants depend on the different dimensions of the cylinder. This flexibility will be needed in Section \ref{secext}.

\begin{lemma}\label{lbnd}
Let $Q$ be a cylinder whose base is given by $B_{R}(x_1)$ in $\left\{t=t_{1}\right\}$ and $B_{R}(x_{2})$ in $\left\{t=t_{2}\right\}$, with $Q\subset Q_{1}$. Let $h=t_{2}-t_{1}>0$, and let $d=|x_{2}-x_{1}|$. Suppose that there exists $\eta, \tau_{1}, \tau_{2}$ so that $\frac{d}{R}\leq \eta$, and $\tau_{1}\leq \frac{h}{R^{2}}\leq\tau_{2}$. Let $u\geq 0$ solve
\begin{equation*}
u_{t}-\MM^{-}(D^{2}u)\geq 0\quad\text{in}\quad Q.
\end{equation*}
Let $\thet\in (0,1)$. Suppose that $u(x, t_{1})\geq 1$ for all $x$ such that $|x-x_{1}|\leq \delta R$. Then there exists $\ga(\theta)$ and $\al=\al(\la, \thet,\eta, \tau_{1}, \tau_{2}, \La, N)>0$,  so that for all $|x-x_{2}|\leq (1-\thet)R$, 
\begin{equation*}
u(x,t_{2})\geq \ga \delta^{\al}.
\end{equation*}
\end{lemma}

\begin{proof}
Without loss of generality, we may perform a transformation to assume that $Q$ is a right cylinder. Suppose we are working on an oblique cylinder. If we let $u^{1}(x,t)=u(x+bt, t)$ (with $b\in\RR^{N}$ to be chosen), then 
\begin{equation}\label{right}
u^{1}_{t}-\MM^{-}(D^{2}u^{1})-b\cdot Du^{1}\geq 0 \quad\text{for}\quad \tilde{Q}= \left\{(x,t): (x+bt, t)\in Q\right\}.
\end{equation}
 
For any oblique cylinder $Q$, we choose $b\in\RR^{N}$ so that $\tilde{Q}$ is a right cylinder. Moreover, we have that the magnitude of $|b|=\frac{d}{h}$.

By scaling and adjusting the operator, we may also assume that $Q=Q_{1}$. We let $u^{2}(x,t)=u^{1}(Rx, ht)$. We have 
\begin{equation*}
\frac{u^{2}_{t}}{h}-\frac{1}{R^{2}}\MM^{-}(D^{2}u^{2})-\frac{1}{R} b\cdot Du^{2}\geq 0.
\end{equation*}
This yields 
\begin{equation}\label{super2}
u^{2}_{t}-\frac{h}{R^{2}}\MM^{-}(D^{2}u^{2})-\frac{h}{R}b\cdot Du^{2}\geq 0 \quad\text{in}\quad Q_{1}.
\end{equation}

Now that we are working in $Q_{1}$, our objective is to obtain an estimate in $B_{(1-\thet)}(0)$ when $t=0$. We fix $x_{0}\in B_{(1-\thet)}(0)$.  Without loss of generality, we may choose $\delta$ so that $\delta\leq \frac{1}{2}\theta$. We examine the cylinder $\tilde{Q}_{\theta, 1}$, which is a cylinder with base $B_{\theta}(0,-1)$ and top $B_{\theta}(x_0, 0)$. By our choices, $\tilde{Q}_{\theta, 1}$ fits inside of $Q_{1}$. Inside of $\tilde{Q}_{\theta, 1}$, $u^{2}$ is a supersolution to \pref{super2}, $u^{2}\geq 0$, and $u^{2}(x,-1)\geq 1$ for all $|x|\leq \delta$. We perform yet another change of coordinates to straighten this cylinder as before. We set $u^{3}(x,t)=u^{2}(x+ct,t)$, where $|c|\leq 1-\thet$. $u^{3}$ solves
\begin{equation}
u^{3}_{t}-\frac{h}{R^{2}}\MM^{-}(D^{2}u^{3})-\left(\frac{h}{R}b+c\right)\cdot Du^{3}\geq 0\quad\text{in}\quad Q_{\theta, 1}=B_{\theta}(0)\times [-1, 0].
\end{equation}
The problem reduces to showing that $u^{3}(0, 0)\geq \ga \delta^{\al}$ for some choice of $\ga, \delta, \al$.

We consider 
\begin{equation}
\psi=\psi(x,t)=((\thet^{2}-\delta^{2})(1+t)-|x|^{2}+\delta^{2})^{2}((\thet^{2}-\delta^{2})(1+t)+\delta^{2})^{-\al}
\end{equation}
where we will choose $\al$ later in the proof. Let 
\begin{equation*}
\hat{Q}=\left\{(x,t): (\thet^{2}-\delta^{2})(1+t)-|x|^{2}+\delta^{2}>0 ~\text{with}~ -1<t<0\right\}. 
\end{equation*}
We note that $\hat{Q}\subset Q_{\thet,1}$, and for $(x,t)\in \partial_{p}\hat{Q} \cap\left\{ t>-1\right\}$ (the lateral boundary), $\psi(x,t)=0$. Moreover, when $t=-1$, we see that $\hat{Q}(t=-1)\subset B_{\delta}(t=-1)$, and $\psi(x, -1)\leq \delta^{4-2\al}$. Therefore, 
\begin{equation*}
\delta^{2\al-4}\psi\leq u^{3} \quad\text{on}\quad \partial_{p}\hat{Q}.
\end{equation*}
Now we are ready to understand the solution properties of $\psi$. We let $\rho=\rho(t)=(\thet^{2}-\delta^{2})(1+t)+\delta^{2}$, and $\vp=\vp(x,t)=\rho(t)-|x|^{2}$, so that $\psi(x,t)=\vp^{2}\rho^{-\al}$. On $\hat{Q}$, after a small calculation, we have
\begin{align*}
&\rho^{\al}\psi_{t}-\frac{h}{R^{2}}\rho^{\al}\MM^{-}(D^{2}\psi)-\rho^{\al}\left(\frac{h}{R}b+c\right)\cdot D\psi \\
&\leq -\frac{\al}{\rho}(\thet^{2}-\delta^{2})\vp^{2}+C_{0}\vp-8\tau_{1}\la |x|^{2}
\end{align*}
where $C_{0}=C_{0}(\theta, \delta, \eta, \tau_{2}, N, \la)=2(\thet^{2}-\delta^{2})+4\left(\eta\theta+(1-\thet)\theta\right)+4\La N\tau_{2}$. 

For certain, if $8\la \tau_{1} |x|^{2}\geq C_{0}\vp$, then $\psi$ is a subsolution with 0 right hand side. If we are in the case where $8\la \tau_{1} |x|^{2}< C_{0}\vp$, then we must have
\begin{align*}
8\la\tau_{1}(\rho-\vp)&<C_{0}\vp\\
\frac{8\la\tau_{1}}{C_{0}+8\la\tau_{1}}\vp&<\rho^{-1}\vp^{2}.
\end{align*}

This yields
\begin{align*}
&\rho^{\al}\psi_{t}-\frac{h}{R^{2}}\rho^{\al}\MM^{-}(D^{2}\psi)-\rho^{\al}\left(\frac{h}{R}b+c\right)\cdot D\psi\\
&\leq \left(-\al(\thet^{2}-\delta^{2})\frac{8\la\tau_{1}}{C_{0}+8\la\tau_{1}}+C_{0}\right)\vp-8\tau_{1}\la |x|^{2}\\
&\leq \left(-\al \frac{6 \thet^{2}\la\tau_{1}}{C_{0}+8\la\tau_{1}}+C_{0}\right)\vp
\end{align*}

using the fact that $\delta\leq \frac{\theta}{2}$. Therefore, if 
\begin{equation}\label{ow}
\al> C_{0}\frac{C_{0}+8\la\tau_{1}}{6 \thet^{2}\la\tau_{1}}\\
\end{equation}
then $\rho^{\al}\psi(x,t)$ is subsolution everywhere in $\hat{Q}$. This is how we will choose $\al$. 

Since $\MM^{-}(\cdot)$ is uniformly elliptic, $\frac{d}{R}$ is bounded, by the comparison principle we must have that
\begin{equation*}
\delta^{2\al-4}\psi\leq u
\end{equation*}
everywhere inside $\hat{Q}$. In particular, we obtain that 
\begin{equation*}
u(0, 0)\geq \delta^{2\al-4}\psi(0,0)=\delta^{2\al-4}\thet^{4}\thet^{-2\al}\geq \theta^{4}\delta^{2\al-4}
\end{equation*}
and this completes the proof. We note that this construction holds for all $x\in B_{1-\theta}(0,0)$ as desired. 
\end{proof}

\begin{remark}\label{rightcyl}
We note that if $\tau_{1}, \tau_{2}, \eta$ are all constants which only depend on universal quantities $\la, \La, N, \theta, \delta$, then there is a lower bound which is universal (this is how the proof is originally presented in \cite{krylov, krylovsaf}. In particular, in the case when $Q=Q_{r}$, an upright cylinder with radius $r$ and height $r^{2}$, then by \pref{ow}, $\al=\al(\la, \La, N, \theta, \delta)$. 
\end{remark}



Next, we state a lemma which will be useful for obtaining local lower bounds when the right hand side is identically 1. 

\begin{lemma}\label{krylov1}
Let $v\in C(\RR^{N+1})$ satisfy
\begin{equation}\label{gf}
\begin{cases}
v_{t}-\MM^{-}(D^{2}v)\geq 1 \quad\text{in}\quad Q_{1},\\
v\geq 0\quad\text{on}\quad \partial_{p}Q_{1}.
\end{cases}
\end{equation}
There exists a constant $C_{k}$ depending on $\ka, N, \la, \La$ so that for all $|x|\leq \ka$,
\begin{equation}
v(x,0)\geq C_{k}. 
\end{equation}
\end{lemma}


\begin{proof}
We compare $v$ to the barrier function $\psi(x,t)=c_{0}(1-|x|^{2}+t)$, with $c_{0}>0$ to be chosen. We point out that $\psi\leq 0$ on $\partial_{p}Q_{1}$ and 
\begin{equation}
\psi_{t}-\MM^{-}(D^{2}\psi)=c_{0}(1+2\La N)\leq 1
\end{equation}
for $c_{0}$ chosen in terms of $\La, N$. By the comparison principle for viscosity solutions, we have that for all $|x|\leq \ka$, 
\begin{equation*}
v(x,0)\geq \psi(x,0)\geq C_{k}
\end{equation*}
where $C_{k}$ depends on $\La, N, \ka$. 
\end{proof}


By scaling, we see that if 
\begin{equation*}
\begin{cases}
v_{t}-\MM^{-}(D^{2}v)\geq \sigma\quad\text{in}\quad Q_{r}(x_{0}, t_{0}), \\
v\geq 0\quad\text{on}\quad \partial_{p}Q_{r}(x_{0}, t_{0}),
\end{cases}
\end{equation*}
then  $v(x,t_{0})\geq C_{k}\sigma r^{2}$ for all $|x-x_{0}|\leq \ka r$. 

Equipped with these results, we are now ready to prove the lower bound for $w(x,t; Q_{r}(x_{0},t_{0}))$. 

\begin{proposition} \label{qlbnd}
Let $Q_{r}(x_{0}, t_{0})\subset Q_{1}$, $\ka\in (0,1)$, and let $w(x,t)$ satisfy
\begin{equation}
\begin{cases}
w_{t}-\MM^{-}(D^{2}w)\geq \chi_{Q_{r}(x_{0}, t_{0})}\quad\text{in}\quad Q_{1}\\
w=0\quad\text{on}\quad \partial_{p}Q_{1}.
\end{cases}
\end{equation}
There exists $C, \rho, \beta>0$ depending only on $\la, \La, N, \ka$ so that for all $|x|\leq \ka$, $0\geq t\geq t_{0}$,
\begin{equation}\label{eqlbnd}
w(x,t)\geq C r^{\rho}\exp(-\beta/r^{2}).
\end{equation}
\end{proposition}

\begin{proof}
We first prove \pref{eqlbnd} in the case when $t=t_{0}$. Let $t^{'}_{0}=t_{0}-\left(\frac{3r^{2}}{4}\right)$. We consider the cylinder $Q_{r/2}(x_{0}, t_{0}^{'})\subset Q_{r}(x_{0}, t_{0})$. By scaling Lemma \ref{krylov1}, since $w\geq 0$, there exists $c_{0}$ so that for all $|x-x_{0}|\leq \frac{r}{4}$,
\begin{equation}\label{started}
w(x,t_{0}^{'})\geq c_{0}r^{2}.
\end{equation} 
Using the information in this disc, we build our way to gaining information in $B_{\ka}(0,t_{0})$. 


If $r>\ka$, we draw an oblique cylinder with base $B_{r}(x_{0}, t_{0}^{'})$, and top $B_{r}(0,t_{0})\supset B_{\ka}(0,t_{0})$. We apply Lemma \ref{lbnd}, with the choices $R=r$, $h\geq \frac{3r^{2}}{4}$, $\eta=\frac{1}{\ka}-1$, $\delta=\frac{1}{4}$, $\theta=\frac{1}{1000}$, and $\tau_{1}=\frac{3}{4}$, and $\tau_{2}=\frac{1}{\ka^{2}}$. There exists a universal constant $C_{e}=C_{e}(\ka, \La, \la, N)$ so that for all $|x|\leq \ka$,
\begin{equation}\label{ez}
w(x,t_{0})\geq C_{e} r^{2},
\end{equation}
and we are done. 

If we are in the case where $r\leq \ka$, then we need to perform an iterative construction. We note that it is enough to show that \pref{eqlbnd} holds for each fixed $y_{0}\in B_{\ka}(0,t_{0})$. Consider the line segment between $(y_{0},t_{0})$ and $(x_{0}, t_{0}^{'})$. Let $\ell$ denote the minimum integer such that $\frac{|x_{0}-y_{0}|}{\ell}\leq \frac{r}{\sqrt{\ell}}$. We note that for $r$ small, since $|x_{0}-y_{0}|\leq 2$, we may choose $\ell\leq 5/r^{2}$. Next, divide the line into $\ell$ segments of equal length, separated by coordinates $(x_{j-1}, t_{j-1})$ and $(x_{j}, t_{j})$.  We may choose $d=|x_{j}-x_{j-1}|\leq \frac{r}{\sqrt{\ell}}$, and $h=|t_{j}-t_{j-1}|=\frac{3r^{2}}{4\ell}$. We then stack a tower of oblique cylinders along this line segment, with base $B_{R}(x_{j}, t_{j})$ and top $B_{R}(x_{j+1}, t_{j+1})$, with $R=\frac{r}{\sqrt{\ell}}$. In each of these cylinders, in the notation of Lemma \ref{lbnd}, $\frac{h}{R^{2}}=\tau_{1}=\tau_{2}$ is constant, and $\frac{d}{R}=\eta$ is constant. Moreover, we may choose $\delta, \theta=\frac{1}{2}.$ Therefore, there exists $c_{1}(\la, \La, N, \ka)<1$ such that for all $|x-x_{j}|\leq \frac{1}{2}R$ 
\begin{equation}
w(x, t_{j})\geq c_{1} w(x, t_{j-1}).
\end{equation}
After $\ell$ iterations,
\begin{equation}
w(y_{0}, t_{0})\geq c_{1}^{\ell} w(x,t_{0}^{'})\geq c_{1}^{\ell} c_{0} r^{2}.
\end{equation}
Relabeling constants as necessary, we see that there exists $\beta=\beta(\la, \La, d,\ka)$ such that for all $y_{0}\in B_{\ka}(0, t_{0})$, 
\begin{equation}\label{easyest}
w(y_{0}, t_{0})\geq Cr^{\rho}\exp(-\beta/r^{2}).
\end{equation}

To prove the estimate for $t\geq t_{0}$, we note at step $\ell-1$, \pref{easyest} yields that for all $|y|\leq \ka-\frac{r}{\sqrt{\ell}}$, 
\begin{equation*}
w\left(y, t_{0}-\frac{r^{2}}{\ell}\right)\geq Cr^{\rho}\exp(-\beta/r^{2}).
\end{equation*}
Using this information, we may obtain the estimate at any $t\geq t_{0}$, $|x|\leq \ka$, by constructing one final, upright, standard cylinder and applying Lemma \ref{lbnd} to \pref{easyest}. We conclude that \pref{easyest} still holds with a constant which differs at most by a constant depending only on $\la, \La, N,\ka$. This completes the proof for all $t\geq t_{0}$.
\end{proof}

\section{Quantitative Lower Bounds for Nonnegative Supersolutions}\label{secgen}
We use the quantitative lower bound on $w(x,t; Q_{r})$ and the Fabes-Stroock estimate to obtain a lower bound on $w(x,t; \left\{f>\al\right\}\cap Q_{r}(x_{0}, t_{0}))$, for $\al>0$. We then compare that to $u$ solving \pref{nonlinear}, with $\al\sim \norm{f}_{L^{N+1}(Q_{1})}$ to obtain Theorem \ref{TSFS}. 

We will refer to the following corollary as the Fabes-Stroock estimate. 
\begin{corollary}\label{classicFS}
Let $E\subset Q_{r}(x_{0}, t_{0})$ such that $\tilde{Q}_{3r}(x_{0}, t_{0})=B_{3r}(x_{0})\times (t_{0}-9r^{2}, t_{0}+9r^{2}]\subset Q_{1}$. For every $\ka\in (0,1)$, there exists $\sigma, C_{cfs}>0$ depending on $\la, \La, N$, such that for all $(x, t)\notin \tilde{Q}_{3r}(x_{0}, t_{0})$, $t\geq t_{0}+9r^{2}$,
\begin{equation}\label{fsineq}
\frac{w(x,t; E)}{w(x,t; Q_{r})}\geq C_{cfs}\left(\frac{|E|}{|Q_{r}|}\right)^{\sigma}.
\end{equation}
\end{corollary}


 Corollary \ref{classicFS} follows from a combination of results. We first consider the linear setting with the operator $\mathcal{L}:=D_{t}-\sum_{i,j} a_{ij}(x, t)D^{2}_{ij}$, where $\left\{a_{ij}(\cdot, \cdot)\right\}$ are uniformly elliptic. By the work of Coiffman and Fefferman \cite[Lemma 5]{coiffef}, a general strategy to prove estimates of the form \pref{fsineq} is to prove a reverse Holder inequality for the Green's function. Indeed, if denote $g(x,y,t,s)$ to be the Green's function corresponding to the operator $\mathcal{L}$, then setting $w(x,t)=w(x,t; E)=\int_{E} g(x,y,t,s)dyds$ solves
\begin{equation*}
\begin{cases}
w_{t}-\sum_{i,j} a_{ij}(x, t)w_{x_{i}x_{j}}=\chi_{E}\quad\text{in}\quad Q_{1},\\
w=0\quad\text{on}\quad \partial_{p}Q_{1}.
\end{cases}
\end{equation*}

However, the result of Coiffman and Fefferman holds if one can show a reverse Holder inequality for all Euclidean balls in a space, whereas the appropriate geometry for us to consider in this problem is cylinders with parabolic scaling. It was shown in a paper of Calderon \cite{cald} that one can extend their methodology to more general metrics, in particular the parabolic metric $\rho((x,t), (y,s))=\sup\left\{\left(\sum_{i}|x_{i}-y_{i}|^{2}\right)^{1/2}, |t-s|\right\}$ whose unit balls are parabolic cylinders. This reduces the proof of Corollary \ref{classicFS} to showing a reverse Holder inequality which holds for parabolic cylinders. 

This reverse Holder inequality was presented by the authors of \cite{amarnorando}:
\begin{theorem}[Amar and Norando, Corollary 2.10, \cite{amarnorando}]\label{esg}
Let $g(x,y,t,s)$ denote the Green's function on $Q_{1}$ corresponding to the operator $D_{t}-\sum_{i,j} a_{ij}(x,t)D^{2}_{ij}$, with $a_{ij}$ smooth in $x,t$. There exists a positive constant $K=K(N, \la, \La)$ such that for every cylinder $Q_{r}(x_{0})\subset \tilde{Q}_{3r}(x_{0}, t_{0})=B_{3r}(x_{0})\times (t_{0}-9r^{2}, t_{0}+9r^{2}]\subset Q_{1}$, we have for all $(x,t)\notin \tilde{Q}_{3r}(x_{0}, t_{0})$, with $t\geq t_{0}+9r^{2}$
\begin{align*}
&\left[ r^{-(N+2)}\int \int_{Q_{r}(x_{0},t_{0})} g(x, t,y,s)^{(N+1)/N} dy ds\right]^{N/(N+1)}\\
&\leq K r^{-(N+2)} \int \int_{Q_{r}(x_{0},t_{0})} g(x, t, y, s) dy ds. 
\end{align*}
\end{theorem}
We point out that the original formulation of \cite{amarnorando} was written for coefficients independent of time. However, in light of the backward boundary Harnack inequality established in \cite{fabessafyuan}, the result easily generalizes to equations with time-dependent coefficients.  Since none of these estimates depend on the smoothness of the coefficients, we may extend them to linear equations with bounded, measurable coefficients by standard approximation arguments (see for example \cite{cabre, cc}). Moreover, once the estimates hold for linear equations with bounded measurable coefficients, they will also hold for equations with Pucci's extremal operators by standard comparison and approximation techniques.  

We will first prove Theorem \ref{LB}, and then show how we may conclude Theorem \ref{TSFS}. 
\begin{proof}[Proof of Theorem \ref{LB}]
We will constantly relabel a constant $c$ when $c$ only depends on dimension. We first prove the estimate for $t=0$. We denote $\tilde{Q}=B_{1-c_{1}m}\times (-1+c_{1}m, -c_{1}m]$. There exists a choice of $c_{1}$, independent of $m$, such that

\begin{equation*}
\frac{|Q_{1}\setminus \tilde{Q}|}{|Q_{1}|}=1-c(1-c_{1}m)^{N} (1-2c_{1}m)\leq \ka m.
\end{equation*}

Therefore,
\begin{equation*}
\frac{|\Ga\cap \tilde{Q}|}{|Q_{1}|}\geq \frac{|\Ga\cap Q_{1}|}{|Q_{1}|}-\frac{|Q_{1}\setminus \tilde{Q}|}{|Q_{1}|}\geq m-\ka m.
\end{equation*}


Next, we cover $\tilde{Q}$ with cylinders $Q_{c_{1}m/4}$, in such way so that each cylinder stays within $Q_{1}$. We claim that there exists at least one smaller cylinder $Q^{*}_{c_{1}m/4}$ such that 

\begin{equation}\label{comcym2}
|\Ga\cap Q^{*}_{c_{1} m/4}| \geq \frac{(m-\ka m}{2} |Q^{*}_{c_{1}m/4}|.
\end{equation}


There exists a covering of $\tilde{Q}$ using less than $2\frac{|Q_{1}|}{|Q_{c_{1}m/4}|}$ small cylinders. If \pref{comcym2} did not hold, then we would have
\begin{align*}
(1-\ka) m \leq \frac{|\Ga\cap \tilde{Q}|}{|Q_{1}|} &\leq \frac{1}{|Q_{1}|} \sum_{Q_{c_{1} m/4}\subset\tilde{Q} } |\Ga\cap Q_{c_{1}m/4}|\\
&\leq \frac{1}{|Q_{1}|}  \frac{2|Q_{1}|}{|Q_{c_{1}m/4}|}\max |\Ga\cap Q_{c_{1} m/4}|\\
&< \frac{1}{|Q_{1}|} \frac{2|Q_{1}|}{|Q_{c_{1}m/4}|} \frac{(1-\ka)m}{2}|Q_{c_{1}m/4}|\\
&< (1-\ka)m
\end{align*}
which is a contradiction. Therefore, we must have that \pref{comcym2} holds in some cylinder $Q^{*}_{c_{1}m/4}$. By construction, $Q^{*}_{c_{1} m/4}$ satisfies the assumptions of Corollary \ref{classicFS}.



By Proposition \ref{qlbnd} and relabeling constants as necessary,
\begin{equation}\label{qclbnd}
w(x,t; Q^{*}_{c_{1}m/4})\geq C [c_{1}m/4]^{\rho}\exp(-\beta/m^{2})\geq Cm^{\rho}\exp(-\beta/m^{2}).
\end{equation}


Letting $E=\Ga\cap Q^{*}_{c_{1}m/4}$ and applying Corollary \ref{classicFS}, \pref{comcym2}, \pref{qclbnd}, and the comparison principle,
\begin{equation}
\frac{u(x,0)}{\al}\geq w(x,0; E) \geq Cm^{\rho}\exp(-\beta/m^{2})[(1-\ka)m]^{\sigma}=cm^{\rho}\exp(-\beta/m^{2})
\end{equation}
for all $|x|\leq \ka$. This establishes the estimate for $t=0$. 

For other times, we have that for any $-\ka m\leq t_{0} \leq 0$, 
\begin{align*}
\frac{|\Ga\cap Q_{1}(0, t_{0})|}{|Q_{1}|}&\geq \frac{|\Ga\cap Q_{1}|}{|Q_{1}|}-\frac{|Q_{1}\setminus (Q_{1}(0, t_{0})\cap Q_{1})|}{|Q_{1}|}\\
&\geq m-\ka m.
\end{align*}

If we extend $w=0$ outside of $Q_{1}$, then $w$ solves \pref{nonlinear} in $Q_{1}(0, t_{0})$, for $E=\Ga\cap Q_{1}(0, t_{0})$. We point out that all of the constants in the estimates above are independent of the choice of $m$. Thus, replacing $m$ by $(1-\ka)m$, we obtain that for all $-\ka m\leq t_{0}\leq 0, |x|\leq \ka$,
\begin{equation}
 u(x, t_{0})\geq w(x,t_{0}; \Ga\cap Q_{1}(0, t_{0}))\geq cm^{\rho}\exp(-\beta/m^{2})\al.
 \end{equation} 
\end{proof}

We now complete the Proof of Theorem \ref{TSFS}.
\begin{proof}[Proof of theorem. \ref{TSFS}]
We note that the right hand side of the estimate is nothing more than the Alexandroff-Backelman-Pucci-Krylov-Tso estimate. For the left hand side, we apply Theorem \ref{LB}. Let $\frac{1}{|Q_{1}|^{1/N+1}}\norm{f}_{L^{N+1}(Q_{1})}=\eta$, then we claim that 
\begin{equation*}
\left|\left\{f> \frac{\eta}{2} \right\}\right|\geq |Q_{1}|\frac{\eta^{N+1}}{2}.
\end{equation*}
To see why, suppose this were not the case. Then we have that since $|Q_{1}|>1$, $f\leq 1$, 
\begin{align*}
\eta^{N+1}=\frac{1}{|Q_{1}|}\int |f|^{N+1}dxdt &\leq \frac{1}{|Q_{1}|} \int |f|^{N+1}dxdt\\
&+\frac{1}{|Q_{1}|}\int_{\left\{f> \eta/2\right\}} |f|^{N+1} dxdt\\
&\leq \left(\frac{\eta}{2}\right)^{N+1}+\frac{\left|{\left\{f> \eta/2\right\}}\right|}{|Q_{1}|}<\eta^{N+1}.
\end{align*}

Therefore, by applying Theorem \ref{LB}, since $\beta\eta^{-(N+1)}\geq \beta|\left\{f>\eta/2\right\}|^{-1}$, and relabeling $\beta$ as necessary, we have
\begin{equation}
u(x,t)\geq c  \left| \left\{ f>\eta/2\right\}\right|^{\rho}\exp(-\beta/|\left\{f>\eta/2\right\}|^{2})\frac{\eta}{2}\geq c \eta^{\rho}\exp(-\beta/\eta^{2(N+1)})
\end{equation}
for all $|x|\leq \ka$, $0\geq t\geq -\ka \frac{\norm{f}^{N+1}_{L^{N+1}(Q_{1})}}{2|Q_{1}|}\geq -\ka \frac{ \left|\left\{f> \frac{\eta}{2}\right\}\right|}{|Q_{1}|}$. 
Relabeling our constants as necessary, this gives us the desired result. 
\end{proof}

\section{Special Cases and Extensions}\label{secext}
In this section, we discuss some useful special cases and extensions of Theorem \ref{LB}. In particular, we obtain estimates comparable to the results of \cite{fabgarsalsa} and \cite{krylovL}, and we recover the estimates of the elliptic case \cite{csw} as $t\rightarrow\infty$. We show several estimates given in terms of $|\left\{f>\al\right\}\cap D|$ for some choice of $D\subset Q_{1}$. 

The first special case we discuss is when $D\subset \left\{t\leq -\ka^{2}\right\}$. 
\begin{corollary}\label{homogest}
Let $u$ be nonnegative in $\overline{S}(f, Q_{1})$. Fix $\ka\in (0,1)$ and set $\Ga=\left\{f>\al\right\}\cap\left\{t\leq -\ka^{2}\right\}$, and $m=\frac{|\Ga|}{|Q_{1}|}$. There exists $c, \rho, \beta>0$ depending only $\la, \La, N, \ka$ so that for all $|x|\leq \ka$, $0\geq t\geq -\ka^{2}/2$,
\begin{equation}
u(x,t)\geq cm^{\rho}\exp((\log 1/m)^{2})\al=cm^{\rho}m^{\beta |\log 1/m|}\al.
\end{equation}
\end{corollary}

This corollary follows in two steps. First, we show that if $Q_{r}(x_{0}, t_{0})\subset \left\{t\leq -\ka^{2}\right\}$, then we can obtain a stronger estimate for all $|x|\leq \ka, t\geq -\ka^{2}/2$. The argument uses a more delicate construction than that which is done in Proposition \ref{qlbnd}. This construction capitalizes on the freedom of parameters of Lemma \ref{lbnd}, and the additional time interval $[-\ka^{2}, -\ka^{2}/2)$ to fit our construction in height-wise. 

\begin{proposition}
Let $Q_{r}(x_{0}, t_{0})\subset (Q_{1}\cap \left\{t\leq -1/2\right\})$, $\ka\in (0,1)$, and let $w$ satisfy
\begin{equation}
\begin{cases}
w_{t}-\MM^{-}(D^{2}w)\geq \chi_{Q_{r}(x_{0}, t_{0})}\quad\text{in}\quad Q_{1},\\
w=0\quad\text{on}\quad \partial_{p}Q_{1}.
\end{cases}
\end{equation}
Then there exists $C, \rho, \beta>0$ depending only on $\la, \La, N, \ka$ so that for all $|x|\leq \ka$, $0\geq t\geq -\frac{\ka^{2}}{2}$,
\begin{equation}
w(x,t)\geq C r^{\rho}\exp(\beta (\log 1/r)^{2}).
\end{equation}
\end{proposition}

\begin{proof}
We first prove the estimate for $t=-\ka^{2}/2$. As in Proposition \ref{qlbnd}, we begin with \pref{started}. If $r>\ka$, then as before there is nothing to prove.  

If we are in the case where $r\leq \ka$, then we perform an iterative construction. We consider a sequence of stacked oblique cylinders with expanding radii, and repeatedly apply Lemma \ref{lbnd} in order to obtain a lower bound for $w$ in $B_{\ka}(0, -\ka^{2}/2)$. Let $\ell$ denote the number of cylinders in our tower. The goal is to choose the final cylinder so that the final radius $R_{\ell}$ satisfies $\frac{1}{2}R_{\ell}>\ka$. As in the construction in Proposition \ref{qlbnd}, given a lower bound in $B_{r/4}$, we can obtain a lower bound at a later time in $B_{r/2}$. Therefore, we choose $R_{j}=2^{j+1}\frac{r}{4}$, and consequently, $\ell=C\log(\ka/r)$. In order to guarantee that the construction fits in height-wise, we may choose $h_{j}=\frac{\ka^{2}}{2\ell}$, for each $j$. Finally, we choose $x_{j}=x_{j-1}\left(1-\frac{R_{j}}{\norm{x_{j-1}}}\right)$ and $t_{j}=t_{j-1}+h_{j}$. We note that there exists a $j^{*}$ such that the line segment between $x_{j^{*}}$ and $x_{j^{*}+1}$ passes through the origin. We then define $x_{k}=0$ for all $k\geq j^{*}$. In the language of Lemma \ref{lbnd}, $\eta<1$, $\tau_{2}=1$, $\tau_{1}\geq \frac{1}{\ell}$ in each iteration, with $\theta=1/2$ and $\delta=1/4$. By \pref{ow}, we see that since $\al\sim \frac{1}{\tau_{1}}$, and for all $|x-x_{j}|\leq R_{j}$, 
\begin{equation}
w(x, t_{j})\geq C\left(\frac{1}{4}\right)^{\ell}w(x, t_{j-1})=Cr^{\beta}w(x, t_{j-1}).
\end{equation}

Therefore, after $\ell$ iterations and relabeling constants as necessary, for all $|x|\leq \ka$, 
\begin{equation}
w(x, -\ka^{2}/2)\geq Cr^{\rho}r^{\beta \log\ka/r}=Cr^{\rho}\exp (\beta (\log 1/r)^{2}).
\end{equation}
where $C, \rho, \beta$ depend on $\la, \La, N, \ka$. A similar argument as in Proposition \ref{qlbnd} yields the estimate for all $t\geq- \frac{\ka^{2}}{2}$. 
\end{proof}

This shows that if $Q_{r}(x_{0}, t_{0})\subset \left\{t\leq -\ka^{2}\right\}$, then the fundamental solution corresponding to $Q_{r}(x_{0},t_{0})$ satisfies a slightly stronger estimate for all $|x|\leq \ka$ and $t\geq -\ka^{2}/2$. Since $\Ga\subset \left\{t\leq -\ka^{2}\right\}$, in the proof of Theorem \ref{LB}, we now have that $Q_{r}(x_{0}, t_{0})\subset \left\{t\leq -\ka^{2}\right\}$. Inserting this into the proof of Theorem \ref{LB}, we obtain that for all $|x|\leq \ka$, $0\geq t\geq -\frac{\ka^{2}}{2}$, 
\begin{equation}
u(x,t)\geq w(x,t)\al\geq Cr^{\rho} r^{\beta |\log (\ka/r)|}\al, 
\end{equation}
which yields Corollary \ref{homogest}. 

Another useful adaptation that we mention is if we are interested in $\Ga=\left\{f>\al\right\}\cap \underline{Q}_{1-\ka}$, where $Q_{1-\ka}=B_{1-\ka}\times \left\{-1\leq t\leq -\ka^{2}\right\}$. This case corresponds to the setting studied in \cite{fabgarsalsa} and \cite{krylovL}, where they assume the set of interest is strictly in the interior of the cylinder. In this case, we are able to obtain a lower bound with a power-type decay. 
\begin{corollary}\label{slicklb}
Let $u$ be nonnegative in $\overline{S}(f, Q_{1})$. Let $\Ga=\left\{f>\al\right\}\cap \underline{Q}_{1-\ka}$, and let $m=\frac{|\Ga|}{|Q_{1}|}$. For every $\ka\in(0,1)$, there exists $c, \beta>0$ which only depend on $\la, \La, N, \ka$ so that for all $|x|\leq \ka$ and $0\geq t \geq -\ka^{2}/2$,
\begin{equation}\label{alg}
u(x,t)\geq cm^{\beta}.
\end{equation}
\end{corollary}


\begin{proof}[Proof of Corollary \ref{slicklb}]
As in the proof of Corollary \ref{homogest}, we only need to show that an estimate of the form \pref{alg} holds if $Q_{r}(x_{0}, t_{0})\subset \underline{Q}_{1-\ka}$. Indeed, we can take $R\geq \ka$ in our construction, $h\geq \ka^{2}/2$, and $\eta\leq \frac{1}{\ka}$. In this way, all of the dependences of Lemma \ref{lbnd} only depend on $\ka$, and ellipticity, and thus after at most two iterations, for all $|x|\leq \ka, 0\geq t\geq -\ka^{2}/2$, 
\begin{equation}
w(x,t)\geq Cr^{2}.
\end{equation}

As in the proof of Corollary \ref{homogest}, this is enough to yield Corollary \ref{slicklb}.
\end{proof}

\begin{remark}
In \cite{fabgarsalsa}, it was shown using the Green's function representation of solutions and the classical covering arguments previously discussed, that for $u\in S(\chi_{\Ga}, B_{r}\times (-\infty, \infty))$ with $r\geq 1$, if $\Ga\subset B_{1/2}\times [0,1]$, then $u(0, 2)\geq c |\Ga|^{m}$ for $c, m$ depending on $\la, \La, N$. From Corollary \ref{slicklb}, we can immediately recover this result. Similarly, in Theorem 4.1 of \cite{krylovL}, the author uses a different approach to show that if $u\in W^{1,2, N+1}(B_{1}\times (0,2))$ is a supersolution to a linear uniformly parabolic nondivergence form equation with smooth coefficients, then there exists $\ga=\ga(\la, \La, d)\in (0,1)$ and $C=C(\la, \La, N)$ such that 
\begin{equation}\label{kryest}
|Q_{1}(0,1)\cap \left\{-Lu\geq \la\right\}|\leq C\la ^{-\ga}u^{\ga}(0,0).
\end{equation}
This can be viewed as an application of Corollary \ref{slicklb} in the cylinder $Q_{2}$ with $\ka=1/2$. 
\end{remark}

We point out that we can also express \pref{bothsides} as
\begin{equation}
c\norm{f}_{L^{N+1}(Q_{1})}^{\al}\leq u(x,t)\leq C \norm{f}_{L^{N+1}(Q_{1})}.
\end{equation}
with $\al=\rho+\beta\norm{f}^{N+1}_{L^{N+1}(Q_{1})}\left|\log \norm{f}_{L^{N+1}(Q_{1})}\right|.$

For comparison, we also state the elliptic version of this result from \cite{csw}:
\begin{theorem}[Caffarelli, Souganidis, Wang, \cite{csw}]\label{ellip}
Let $u\geq 0$ solve $-\MM^{-}(D^{2}u) \\ \geq f$ in $B_{1}$, $u=0$ on $\partial B_{1}$ with $0\leq f\leq 1$. For every $\ka\in(0,1)$, there exists $c, C,\al>0$ depending only on $\la, \La, N, \ka$ such that for all $|x|\leq \ka$, 
\begin{equation*}
c\norm{f}_{L^{N}(B_{1})}^{\al}\leq u(x)\leq C\norm{f}_{L^{N}(B_{1})}.
\end{equation*}
\end{theorem}

\begin{remark}
We may obtain Theorem \ref{ellip} by looking at Corollary \ref{homogest}. In the time-independent setting, we consider the limit as $t\rightarrow \infty$, and a cylinder of infinite height with $f$ constant in time.  Since we do not worry about fitting our construction in heightwise, we may always choose cylinders in our tower with dimensions so that $\tau_{1}, \tau_{2}$ and $\eta$ depend only on $\la, \La, N, \ka$.  Therefore, for all $|x|\leq \ka$, $w(x,t)\geq c_{1}^{\ell}c_{0}r^{2}=cr^{\beta}$, where $\beta=\beta(\la, \La, N, \ka)$. This is enough to conclude. 
\end{remark}

\begin{remark}
We mention that an alternative proof of Nikolai Krylov \cite{krylovjess}, does not require the Fabes-Stroock estimate, and yields a stronger lower bound in the general setting. 
\end{remark}

\section*{Acknowledgements} This work was completed as a part of the author's doctoral thesis. The author would like to thank her thesis advisor, Takis Souganidis, for his patient guidance and many helpful discussions. Also, the author would like to thank Carlos Kenig, for referring her to \cite{cald} in order to obtain the Fabes-Stroock estimate for parabolic cylinders. The author would finally like to thank Nikolai Krylov, for suggesting an improvement to Proposition \ref{qlbnd}, and for teaching her techniques and results from parabolic regularity theory. The author was supported by the NSF Graduate Research Fellowship on NSF grant DGE-1144082.

\bibliographystyle{amsplain}
\bibliography{regnodea}

\end{document}